\documentclass[11pt]{amsart}

\usepackage{amsmath}
\usepackage{amssymb}

\newtheorem{theorem}{Theorem}[section]
\newtheorem{claim}[theorem]{Claim}

\newtheorem{sclaim}[theorem]{SubClaim}

\theoremstyle{definition}
\newtheorem{definition}[theorem]{Definition}

\theoremstyle{remark}
\newtheorem{remark}[theorem]{Remark}

\newtheorem{conclusion}[theorem]{Conclusion}

\newtheorem{notation}[theorem]{Notation}

\newcount\skewfactor
\def\mathunderaccent#1#2 {\let\theaccent#1\skewfactor#2
\mathpalette\putaccentunder}
\def\putaccentunder#1#2{\oalign{$#1#2$\crcr\hidewidth
\vbox to.2ex{\hbox{$#1\skew\skewfactor\theaccent{}$}\vss}\hidewidth}}

\def\smallbox#1{\leavevmode\thinspace\hbox{\vrule\vtop{\vbox
   {\hrule\kern1pt\hbox{\vphantom{\tt/}\thinspace{\tt#1}\thinspace}}
   \kern1pt\hrule}\vrule}\thinspace}

\newcommand{\bool}{{\bf B}}

\newcommand{\cf}{{\rm cf}}

\newcommand{\Depth}{{\rm Depth}}

\newcommand{\st}{{such that}}
\newcommand{\seq}{{sequence}}
\newcommand{\cont}{{continuous}}
\newcommand{\incr}{{increasing}}
\newcommand{\Then}{{\underline{Then}}}

\newcommand{\mat}{\mathcal}

\def\qedref#1{$\qed_{\reforiginal{#1}}$\marginpar{$\langle$#1$\rangle$}}

\title{On ${\rm DEPTH}$ and ${\rm DEPTH}^+$ of Boolean Algebras}
\author{Shimon Garti and Saharon Shelah}
\address{Institute of Mathematics
 The Hebrew University of Jerusalem
 Jerusalem 91904, Israel
 and  Department of Mathematics
 Rutgers University
 New Brunswick, NJ 08854, USA}
\email{shelah@math.huji.ac.il}
\urladdr{http://www.math.rutgers.edu/\char`\~shelah}
\thanks{First typed: August 2000 \\
Research supported by the United States-Israel Binational
Science Foundation. Publication 878}
\subjclass{}
\keywords{Boolean Algebras, Depth}

\begin{document}

\let\labeloriginal\label
\let\reforiginal\ref

\begin{abstract}
%%% Put abstract here:
We show that the ${\rm Depth}^+$ of an ultraproduct
of Boolean Algebras, can not jump over the ${\rm Depth}^+$
of every component by more than one cardinality. 
We can have, consequently, similar results for the Depth 
invariant.
\end{abstract}

\maketitle

% start document here:
%\newpage

\section{introduction}
Monk \cite{Mo96} has dealt systematically with 
cardinal invariants of Boolean algebras.
In particular he dealt with \smallbox{the question}
how an invariant of an unltraporduct of a sequence 
of Boolean algebras relate to the ultraproduct of the sequence of 
the invariants of each of the Boolean algebras.
That is the relationship of 
${\rm inv} (\prod_{\epsilon<\kappa}
\bool_\epsilon \setminus D)$ with $\prod_{\epsilon<\kappa}
{\rm inv} (\bool_\epsilon)\setminus D$. One 
of the invariants he dealt with is the depth of a 
Boolean algebra, ${\rm Depth} (\bool)$. We continue 
here \cite{Sh:853} getting weaker results without ``large 
cardinal axioms''. On related results see 
\cite{MgSh:433}, \cite{Sh:703}, \cite{RoSh:733}.
Further results on Depth and ${\rm Depth}^+$ 
by the authors are in preparation are
in \cite{Sh:F754}.

Recall: 
\begin{definition}
\label{0.1}
Let ${\bf B}$ be a Boolean Algebra.
$$
\Depth(\bool):= \sup\{\theta:\exists \bar b=(b_\gamma:\gamma<\theta), 
\hbox{ \incr\ \seq\ in } \bool\}
$$
Dealing with questions of Depth, Saharon Shelah noticed that investigating 
a slight little modification of Depth, namely - $\Depth^+$, might be helpful 
(see \cite{Sh:853} for the behavior of Depth and $\Depth^+$ above a compact cardinal).
\end{definition}

Recall:

\begin{definition}
\label{0.2}
Let $\bool$ be a Boolean Algebra.
$$
\Depth^+ (\bool):=\sup \{\theta^+:\exists \bar b=(b_\gamma:\gamma<\theta), 
\hbox{ \incr\ \seq\ in } \bool\}
$$
This article deals mainly with $\Depth^+$, in the aim to get results for the Depth.
It follows \cite{Sh:853}, both - in the 
general ideas and in the method of the proof.

Let us take a look on the main claim of \cite{Sh:853}:
\end{definition}

\begin{claim}
\label{0.3}
Assume
\begin{enumerate}
\item[(a)] $\kappa<\mu\leq \lambda$
\item [(b)] $\mu$ is a compact cardinal
\item[(c)] $\lambda=\cf (\lambda)$
\item[(d)] $(\forall \alpha<\lambda) (|\alpha|^\kappa<\lambda)$
\item[(e)] $\Depth^+ (\bool_i)\leq \lambda$, for every $i<\kappa$
\end{enumerate}
\Then\ $\Depth^+ (\bool)\leq \lambda$

{\rm So, $\lambda$ bounds the $\Depth^+ (\bool)$, 
where $\bool$ is the ultraproduct of 
the Boolean Algebras $\bool_i$, 
if it bounds the $\Depth^+$ of every $\bool_i$. 
That requires some reasonable assumptions on $\lambda$, 
and also a pretty high price for 
that result - you should raise your view to 
a very large $\lambda$, above a compact cardinal. 
Now, the existence of large cardinals 
is an interesting philosophical question. 
You might think that adding a compact cardinal 
to your world is a natural extending of ZFC.
But, mathematically, it is important to check 
what happens without a compact cardinal 
(or bellow the compact, even if the compact cardinal exists).

In this article we drop the assumption of compact cardinal. 
Consequently, we phrase a weaker conclusion. We prove that if 
$\lambda$ bounds the $\Depth^+$ of every 
$\bool_i$, then the $\Depth^+$ of $\bool$
can not jump beyond $\lambda^+$.}
\end{claim}
%\newpage

\section{Bounding ${\rm DEPTH}^+$}

\begin{notation}
\label{1.1}
\begin{enumerate}
\item[(a)] $\kappa,\lambda$ are infinite cardinals
\item[(b)] $D$ is an ultrafilter on $\kappa$
\item[(c)] $\bool_i$ is a Boolean Algebra, for any $i<\kappa$
\item[(d)] $\bool=\prod\limits_{i<\kappa} \bool_i /D$
\end{enumerate}
\end{notation}

\begin{claim}
\label{1.2}
Assume 
\begin{enumerate}
\item[(a)] $\lambda=\cf(\lambda)$
\item[(b)] $(\forall \alpha<\lambda) (|\alpha|^\kappa<\lambda)$
\item[(c)] $\Depth^+ (\bool_i) \leq \lambda$, for every $i<\kappa$
\end{enumerate}
\Then\ $\Depth^+ (\bool) \leq \lambda^+$
\end{claim}

\begin{conclusion}
\label{1.3}
Assume 
\begin{enumerate}
\item[(a)] $\lambda^\kappa=\lambda$
\item[(b)] $\Depth (\bool_i) \leq \lambda$, for every $i<\kappa$
\end{enumerate}
\Then\ $\Depth (\bool) \leq \lambda^+$
\end{conclusion}

\begin{proof} 
By part (b), $\Depth^+ (\bool_i)\leq \lambda^+$ for every $i<\kappa$.
By part (a), $\alpha<\lambda^+ \Rightarrow |\alpha|^\kappa<\lambda^+$. 
Now, $\lambda^+$ is a regular cardinal, so the pair $(\kappa,\lambda^+)$ 
satisfies the requirements of claim \ref{1.2}. 
So, $\Depth^+ (\bool)\leq \lambda^{+2}$, and that means that $\Depth (\bool) \leq
\lambda^+$.
\hfill \qedref{1.3}
\end{proof}

\begin{remark}
\label{1.4}
If $\lambda$ is inaccessible 
(or even strong limit, with cofinality above 
$\kappa$), and $\Depth (\bool_i)<\lambda$ for every 
$i<\kappa$, you can verify easily that $\Depth (\bool)<\lambda$, by claim 
\ref{1.2} and simple cardinal arithmetic.
\end{remark}

\par \noindent 
Proof of claim \ref{1.2}:\ 
Let $\langle M_\alpha:\alpha<\lambda^+\rangle$ be continuous and increasing \seq\ 
of elementary submodels of $({\mat H} (\chi),\in)$ for large enough 
$\chi$, with the next properties: 
\begin{enumerate}
\item[(a)] $(\forall \alpha <\lambda^+) (\| M_\alpha\|=\lambda)$
\item[(b)] $(\forall \alpha<\lambda^+) (\lambda+1
\subseteq M_\alpha)$
\item[(c)] $(\forall \beta<\lambda^+) (\langle M_\alpha:\alpha
\leq \beta\rangle \in M_{\beta+1})$
\end{enumerate}
Choose $\delta^*\in S^{\lambda^+}_\lambda (:= \{\delta<\lambda^+:
\cf (\delta)=\lambda)$, \st\ $\delta^*=M_{\delta^*} \cap \lambda^+$.
Assume toward contradiction, that $(a_\alpha:\alpha<\lambda^+)$ is an
\incr\ \seq\ in $\bool$. Let us write $a_\alpha$ as 
$\langle a^\alpha_i:i<\kappa\rangle/D$ for every $\alpha<\lambda^+$.
We may assume that $\langle a^\alpha_i:
\alpha<\lambda^+, i<\kappa\rangle \in M_0$.

We will try to create a set $Z$, in the subclaim bellow, with 
the following properties: 
\begin{enumerate}
\item[(a)] $Z \subseteq \lambda^+, |Z|=\lambda$
\item[(b)] $\exists i_* \in \kappa$ \st\ for every $\alpha<\beta,
\alpha,\beta\in Z$, 
we have $\bool_{i_*} \models a^\alpha_{i_*} 
<a^\beta_{i_*}$
\end{enumerate}
Since $|Z|=\lambda$, we have an \incr\ \seq\ of length $\lambda$
in $\bool_{i_*}$, so $\Depth^+ (\bool_{i_*}) \geq \lambda^+$, contradicting 
the assumptions of the claim.
\medskip
\hfill \qedref{1.2}

\begin{sclaim}
\label{1.5}
There exists $Z$ as above
\end{sclaim}

\begin{proof}
For every $\alpha<\beta<\lambda^+$, define:
$$
A_{\alpha,\beta}=\{i<\kappa:\bool_i \models a^\alpha_i<a^\beta_i\}
$$
By the assumption, $A_{\alpha,\beta} \in D$, 
for any $\alpha<\beta<\lambda^+$. 
For any $\alpha<\delta^*$
Let $A_\alpha$ denote the set $A_{\alpha,\delta^*}$. \newline
Let 
 $\langle v_\alpha:\alpha<\lambda\rangle$ be \incr\ and \cont, \st:
\begin{enumerate}
\item[(i)] $v_\alpha\in [\delta^*]^{<\lambda}$, for every 
$\alpha<\lambda$
\item[(ii)] $v_\alpha$ has no last element, for every 
$\alpha<\lambda$
\item[(iii)] $\delta^*=\bigcup\limits_{\alpha<\lambda} 
v_\alpha$
\end{enumerate}
Let $u\subseteq \delta^*$, 
$|u|\leq \kappa$. Define:
$$
S_u=\{\beta<\delta^*:\beta>{\rm sup} (u) \hbox{ and } 
(\forall \alpha\in u) (A_{\alpha,\beta}=A_\alpha)\}.
$$ 
Now, define $C=\{\delta<\lambda:\delta$ 
is a limit ordinal, and 
$$
(\forall \alpha<\delta)[(u\subseteq v_\alpha)\wedge 
(|u|\leq \kappa)\Rightarrow {\rm sup}(v_\delta)={\rm sup}
(S_u\cap {\rm sup} (v_\delta))]\}
$$
Since $\lambda=\cf(\lambda)$ and $(\forall \alpha<\lambda) 
(|\alpha|^\kappa <\lambda)$, and since $|v_\delta|<\lambda$, 
clearly $C$ is a club set of $\lambda$. 

The fact that $|D|=2^\kappa <\cf(\lambda)=\lambda$ 
implies that there exists $A_*\in D$ \st\ 
$S=\{\alpha<\lambda:\cf (\alpha)>\kappa$ and 
$A_{{\rm sup}(v_\alpha)}=A_*\}$ is a stationary 
subset of $\lambda$.

$C$ is a club and $S$ is stationary, so $C\cap S$ is also 
stationary. Choose $\delta^1_0=0$. 
Choose $\delta^1_{\epsilon+1}\in 
C\cap S$ for every $\epsilon<\lambda$, \st\ 
$\epsilon<\zeta\Rightarrow {\rm sup}
\{\delta^1_{\epsilon+1}:\epsilon<\zeta\}<\delta^1_{\zeta+1}$. 
Define $\delta^1_\epsilon$ to be the 
limit of $\delta^1_{\gamma+1}$, when $\gamma<\epsilon$, 
for every limit $\epsilon<\lambda$. Since $C$ is closed, 
we have:
\begin{enumerate}
\item[(a)] $\{\delta^1_\epsilon:\epsilon<\lambda\}
\subseteq C$
\item[(b)] $\langle \delta^1_\epsilon:\epsilon<\lambda\rangle$ 
is \incr\ and \cont\ 
\item[(c)] $\delta^1_{\epsilon+1} \in S$, for every $\epsilon<\lambda$
\end{enumerate}
Lastly, define $\delta^2_\epsilon={\rm sup} 
(v_{\delta^1_\epsilon})$, for every 
$\epsilon<\lambda$. Define, for every $\epsilon<\lambda$, the next family:
$$
{\frak A}_\epsilon=\{S_u\cap \delta^2_{\epsilon+1}\setminus 
\delta^2_\epsilon:u\in [v_{\delta^2_{\epsilon+1}}]^{\leq \kappa}\}
$$
We get a family of non-empty sets, which is downward $\kappa^+$-directed.
So, there is a $\kappa^+$-complete filter $E_\epsilon$ on
$[\delta^2_\epsilon, \delta^2_{\epsilon+1})$, with ${\frak A}_\epsilon
\subseteq E_\epsilon$, for every $\epsilon<\lambda$.

Define, for any $i<\kappa$ and $\epsilon<\lambda$, the sets
$W_{\epsilon,i} \subseteq [\delta^2_\epsilon,\delta^2_{\epsilon+1})$ 
and $B_\epsilon\subseteq \kappa$, by:
$$
W_{\epsilon,i} := \{\beta:\delta^2_\epsilon\leq \beta<\delta^2_{\epsilon+1} 
\hbox{ and } i\in A_{\beta, \delta^2_{\epsilon+1}}\}
$$
$$
B_\epsilon:=\{i<\kappa:W_{\epsilon,i}\in E^+_\epsilon\}
$$
At last, take a look on $W_\epsilon:=\cap \{[\delta^2_\epsilon,
\delta^2_{\epsilon+1})\setminus W_{\epsilon,i}:
i\in \kappa \setminus B_\epsilon\}$.
For every $\epsilon<\lambda, W_\epsilon\in E_\epsilon$, since 
$E_\epsilon$ is $\kappa^+$-complete, so clearly $W_\epsilon\neq \emptyset$.

Choose 
$\beta=\beta_\epsilon\in W_\epsilon$. 
If $i\in A_{\beta,\delta^2_{\epsilon+1}}$, 
then $W_{\epsilon,i} \in E^+_\epsilon$, so $A_{\beta,\delta^2_{\epsilon+1}}\subseteq 
B_\epsilon$ (by the definition of $B_\epsilon$). But, 
$A_{\beta,\delta^2_{\epsilon+1}}\in D$, so 
$B_\epsilon\in D$, and consequently - $A_*\cap B_\epsilon\in D$, for any 
$\epsilon<\lambda$. \newline
Choose $i_\epsilon\in A_*\cap B_\epsilon$, for every $\epsilon<\lambda$. 
You choose $\lambda\  i_\epsilon$-s from $A_*$, and $|A_*|=\kappa$, so we 
can arrange a fixed $i_*\in A_*$ \st\ the set $Y=\{\epsilon<\lambda:\epsilon$ 
is even ordinal, and $i_\epsilon=i_*\}$ has cardinality $\lambda$.

The last step will be as follows: \newline
define $Z=\{\delta^2_{\epsilon+1}:\epsilon\in Y\}$. 
Clearly, $Z\in [\delta^*]^\lambda \subseteq [\lambda^+]^\lambda$.
We will show that for $\alpha<\beta$ from $Z$, we get 
$\bool_{i_*} \models a^\alpha_{i_*}<a^\beta_{i_*}$. 
The idea is that if $\alpha<\beta$ and $\alpha,\beta\in Z$,
then $i_*\in A_{\alpha,\beta}$.

Why? Recall that $\alpha=\delta^2_{\epsilon+1}$ and 
$\beta=\delta^2_{\zeta+1}$, for some $\epsilon<\zeta<\lambda$ 
(that's the form of the members of $Z$). Define: \newline
$U_1=S_{\{\delta^2_{\epsilon+1}\}} \cap 
[\delta^2_\zeta,\delta^2_{\zeta+1}) \in 
{\frak A}_\zeta \subseteq E_\zeta$ \newline
$U_2=\{\gamma:\delta^2_\zeta\leq \gamma<\delta^2_{\zeta+1}$ and 
$i_*\in A_{\gamma,\delta^2_{\zeta+1}}\}\in E^+_\zeta$ \newline
So, $U_1\cap U_2\neq \emptyset$. \newline
Choose, $\iota\in U_1\cap U_2$. 
\begin{enumerate}
\item[(a)] $\bool_{i_*} \models a^\alpha_{i_*}<a^\iota_{i_*}$

[Why? well, $\iota\in U_1$, so $A_{\delta^2_{\epsilon+1,\iota}}=
A_{\delta^2_{\epsilon+1}}=A_*$. But, $i_*\in A_*$, so 
$i_*\in A_{\delta^2_{\epsilon+1,\iota}}$, which means that 
$\bool_{i_*} \models a^{\delta^2_{\epsilon+1}}_{i_*}
 (=a^\alpha_{i_*})<a^\iota_{i_*}$].
\item[(b)] $\bool_{i_*}\models a^\iota_{i_*}<a^\beta_{i_*}$

[Why? well, $\iota\in U_2$, so $i_*\in A_{\iota,\delta^2_{\zeta+1}}$, 
which means that $\bool_{i_*} \models a^\iota_{i_*}< 
a^{\delta^2_{\zeta+1}}_{i_*} (=a^\beta_{i_*})$].
\item[(c)] $\bool_{i_*}\models a^\alpha_{i_*}<a^\beta_{i_*}$

[Why? by (a)+(b)].
\end{enumerate}
So, we are done. 
\end{proof}
\hfill \qedref{1.5}

Without a compact cardinal, we may have a `jump' of the $\Depth^+$ 
in the ultraproduct of the Boolean Algebras (see \cite[\S5]{Sh:652}). 
So, we can have $\kappa<\lambda$,
$\Depth^+ (\bool_i)\leq \lambda$ for every $i<\kappa$, and 
$\Depth^+(\bool)=\lambda^+$. We can show that if there exists such an 
example for $\kappa$ and $\lambda$, then you can create an example 
for every regular $\theta$ between $\kappa$ and $\lambda$.

\begin{claim}
\label{1.6} 
Assume 
\begin{enumerate}
\item[(a)] $\kappa<\lambda, D$ is an ultrafilter on $\kappa$
\item[(b)] $\Depth^+ (\bool_i)\leq \lambda$, for every $i<\kappa$
\item[(c)] $\Depth^+ (\bool)=\lambda^+$
\item[(d)] $\theta\in {\rm Reg} \cap [\kappa,\lambda)$
\end{enumerate}
\Then\ there exists Boolean Algebras ${\bf C}_j$ for any 
$j<\theta$, and a uniform 
ultrafilter $E$ on $\theta$, \st\ $\Depth^+
({\bf C}_j)\leq \lambda$ for every $j<\theta$, and 
$\Depth^+ ({\bf C}):=\Depth^+ (\prod\limits_{j<\theta} 
{\bf C}_j / E)=\lambda^+$.
\end{claim}

\begin{proof}
Break $\theta$ into $\theta$ sets 
$(u_\alpha:\alpha<\theta)$, \st:
\begin{enumerate}
\item[(a)] $|u_\alpha|=\kappa$, for every $\alpha<\theta$
\item[(b)] $\bigcup\limits_{\alpha<\theta} u_\alpha=\theta$
\item[(c)] $\alpha\neq \beta\Rightarrow u_\alpha\cap 
u_\beta=\emptyset$
\end{enumerate}
For every $\alpha<\theta$, let $f_\alpha:\kappa \rightarrow 
u_\alpha$ be one to one, onto and order preserving. 
Define $D_\alpha$ on $u_\alpha$, in the following way: 
If $A \subseteq u_\alpha$ then $A \in D_\alpha$ iff 
$f^{-1}_\alpha (A)\in D$. For $\theta$ itself, define a filter 
$E_*$ on $\theta$ in the next way: 
If $A \subseteq \theta$, then $A\in E_*$ iff $A \cap 
u_\alpha \in D_\alpha$ for every (except, maybe 
$<\theta$ ordinals) $\alpha<\theta$. 
Now, choose any ultrafilter $E$ on $\theta$, \st\ $E_*
\subseteq E$.

Define ${\bf C}_{f_\alpha(i)} = \bool_i$, for every 
$\alpha<\theta$ and $i<\kappa$. You will get 
$({\bf C}_j:j<\theta)$ \st\ $\Depth^+ ({\bf C}_j)\leq \lambda$ 
for every $j<\theta$. But, we will show that $\Depth^+ ({\bf C})
\geq \lambda^+$ (remember that ${\bf C}=\prod\limits_{j<\theta}
{\bf C}_j/E$).

Well, let $(a_\xi:\xi<\lambda)$ testify $\Depth^+ (\bool)
=\lambda^+$. Recall, $a_\xi$ is $\langle a^\xi_i:i<\kappa\rangle
/D$. We may write $f_\alpha (a_\xi)$ for $\langle f_\alpha 
(a^\xi_i):i<\kappa\rangle /D_\alpha$, when 
$\alpha<\theta$. Clearly, $(f_\alpha (a_\xi):\xi<\lambda)$ 
testify $\Depth^+ ({\bf C}^\alpha)=\lambda^+$, when 
${\bf C}^\alpha := \prod\limits_{i<\kappa} {\bf C}_{f_\alpha(i)}
/D_\alpha$.

Now, $\langle (f_\alpha (a_\xi):\alpha<\theta):\xi<\lambda\rangle/E$ 
is an \incr\ \seq\ in ${\bf C}$.
\end{proof}
\hfill \qedref{1.6}

\begin{remark}
\label{1.7}
\begin{enumerate}
\item Claim \ref{1.6} applies, in a similar fashion, to the Depth 
invariant.
\item Claim \ref{1.6} is meaningful for comparing $\Depth ({\bf C})$ 
to $\prod\limits_{j<\theta} \Depth ({\bf C}_j)/E$, when 
$\lambda^\theta=\lambda$.
\end{enumerate}
\end{remark}

\end{document}